\RequirePackage{fix-cm}
\documentclass[envcountsect]{svjour3}                     

\smartqed  
\usepackage{graphicx}
\usepackage{caption}
\DeclareCaptionLabelSeparator{twospace}{\ ~}
\usepackage{amsfonts}
\usepackage{mathptmx}
\usepackage[]{geometry}
\spnewtheorem{fact}{Fact}[section]{\bf}{\it}

\begin{document}
\title{On the largest graph-Lagrangians of  hypergraphs
}


\author{Qingsong Tang         \and   Yuejian  Peng  \and   Xiangde Zhang \and Cheng Zhao
}


\institute{Qingsong Tang \at
              College of Sciences, Northeastern University, Shenyang, 110819, P.R.China. \at
               School of Mathematics, Jilin University, Changchun, 130012, P.R.China.\\
              \email{t\_qsong@sina.com}
              \and
              Yuejian  Peng \at
              College of Mathematics, Hunan University, Changsha 410082, P.R. China. \at This research is supported by National Natural Science Foundation of China (No. 11271116).\\
              \email{ypeng1@163.com}
               \and
           Xiangde Zhang (Corresponding author) \at
             College of Sciences, Northeastern University, Shenyang, 110819, P.R.China\\
             \email{zhangxdneu@163.com}
              \and
             Cheng Zhao\at
             Department of Mathematics and Computer Science, Indiana State University, Terre Haute, IN, 47809 USA. \at School of Mathematics, Jilin University, Changchun 130012, P.R. China.\\
             \email{cheng.zhao@indstate.edu}
}

\date{Received: date / Accepted: date}

\maketitle

\begin{abstract}
Frankl and F\"{u}redi (1989) conjectured that the $r$-graph with $m$ edges formed by taking the first $m$ sets in the colex ordering of ${\mathbb N}^{(r)}$ has the largest graph-Lagrangian of all $r$-graphs with $m$ edges.  In this paper, we establish some bounds   for graph-Lagrangians of some special $r$-graphs that support this conjecture.

\end{abstract}

\section{Introduction}
\label{intro}
A useful tool in extremal problems of   hypergraphs (graphs) is the graph-Lagrangian of a   hypergraph (graph) defined as follows.
\begin{definition}
For  an $r$-uniform graph $G$ with the vertex set $[n]$, edge set $E(G)$, and a vector $\vec{x}=(x_1,\ldots,x_n) \in {\mathbb R}^n$,
we associate a homogeneous polynomial in $n$ variables, denoted by  $\lambda (G,\vec{x})$  as follows:
$$\lambda (G,\vec{x})=\sum_{i_1i_2 \cdots i_r \in E(G)}x_{i_1}x_{i_2}\ldots x_{i_r}.$$
Let $S=\{\vec{x}=(x_1,x_2,\ldots ,x_n): \sum_{i=1}^{n} x_i =1, x_i
\ge 0 {\rm \ for \ } i=1,2,\ldots , n \}$.
Let $\lambda (G)$ represent the maximum
 of the above homogeneous  multilinear polynomial of degree $r$ over the standard simplex $S$. Precisely
 $$\lambda (G)= \max \{\lambda (G, \vec{x}): \vec{x} \in S \}.$$
\end{definition}
The value $x_i$ is called the {\em weight} of the vertex $i$.
A vector $\vec{x}=(x_1, x_2, \ldots, x_n) \in {\mathbb R}^n$ is called a feasible weighting for $G$ if
$\vec{x}\in S$. A vector $\vec{y}\in S$ is called an {\em optimal weighting} for $G$
if $\lambda (G, \vec{y})=\lambda(G)$. We call $\lambda(G)$   the  graph-Lagrangian of $G$.



\begin{remark}
$\lambda (G)$ was called Lagrangian of $G$ in literature  (Peng and Zhao 2013; Talbot 2002; Frankl and R\"{o}dl 1984; Frankl and F\"{u}redi 1989, etc.). The terminology `graph-Lagrangian' was suggested by Franco Giannessi.
\end{remark}
 The following fact is easily implied by the definition of the graph-Lagrangian.

\begin{fact}\label{mono}
Let $G_1$, $G_2$ be $r$-graphs and $G_1\subseteq G_2$. Then $\lambda (G_1) \le \lambda (G_2).$
\end{fact}

Graph-Lagrangians of graphs and hypergraphs is one of the most fundamental tools in extremal combinatorics.  Applications of graph-Lagrangian method can be found in  (Mubayi 2006; Frankl and R\"{o}dl 1984; Keevash 2011; Frankl and F¨¹redi 1989).   One of the most intriguing and most simply formulated problems in this topic is the Frankl-F\"{u}redi conjecture that asks how large can the graph-Lagrangian of an $r$-graph with $m$ edges be?
For distinct $A, B \in {\mathbb N}^{(r)}$, we say that $A$ is less than $B$ in the {\em colex ordering} if $max(A \triangle B) \in B$, where $A \triangle B=(A \setminus B)\cup (B \setminus A)$. For example, we have $246 < 150$ in ${\mathbb N}^{(3)}$ since $max(\{2,4,6\} \triangle \{1,5,6\}) \in \{1,5,6\}$.  Note that the first $t \choose r$ $r$-tuples in the colex ordering of ${\mathbb N}^{(r)}$ are the edges of $[t]^{(r)}$. The following conjecture of Frankl and F\"uredi (if it is true) proposes a  solution to the question mentioned at the beginning.

\begin{conjecture} (Frankl and F\"{u}redi 1989)\label{conjecture} The $r$-graph with $m$ edges formed by taking the first $m$ sets in the colex ordering of ${\mathbb N}^{(r)}$ has the largest graph-Lagrangian of all $r$-graphs with  $m$ edges. In particular, the $r$-graph with $t \choose r$ edges and the largest graph-Lagrangian is $[t]^{(r)}$.
\end{conjecture}

We know very little about this conjecture for $r\geq 3$.
It has not been solved even for the case $r = 3$ despite 25 years of intensive research and despite the fact that the problem is coming with quite a concrete and natural conjecture as to the structure of extremal $r$-graphs. There have been, however, some partial results.
Let $C_{r,m}$ denote the $r$-graph with $m$ edges formed by taking the first $m$ sets in the colex ordering of ${\mathbb N}^{(r)}$. The following result was given by  Talbot (2002).
\begin{lemma} (Talbot 2002) \label{LemmaTal7}
For any integers $m, t,$ and $r$ satisfying ${t-1 \choose r} \le m \le {t-1 \choose r} + {t-2 \choose r-1}$,
we have $\lambda(C_{r,m}) = \lambda\left([t-1]^{(r)}\right)$.
\end{lemma}
For the case $r=3$, Talbot(2002) proved the following.
\begin{theorem}  (Talbot 2002)  \label{Talbot}Let $m$ and $t$ be integers satisfying
${t-1 \choose 3} \le m \le {t-1 \choose 3} + {t-2 \choose 2} - (t-1).$
Then Conjecture \ref{conjecture} is true for $r=3$ and this value of $m$.
Conjecture \ref{conjecture} is also true for $r=3$ and $m= {t \choose 3}-1$ or $m={t \choose 3} -2$.
\end{theorem}
For general  $r \ge 4$, the following result is  also proved in (Talbot 2002), which is the evidence for Conjecture \ref{conjecture}  for $r$-graphs $G$ on exactly $t$ vertices.
\begin{theorem} (Talbot 2002)  \label{Talr} For any integer $r \ge 4$ there exists constants $\gamma_r$ and $\kappa_0(r)$ such that if $m$ satisfies
$${t-1 \choose r} \le m \le {t-1 \choose r} + {t-2 \choose r-1} - \gamma_r (t-1)^{r-2}$$
with $t \ge \kappa_0(r)$, let $G$ be an $r$-graph on $t$ vertices with $m$ edges, then $\lambda(G)\leq\lambda([t-1]^{(r)})$.
\end{theorem}
Let us introduce a useful definition  called left-compressed.
\begin{definition}(Talbot 2002)
An $r$-graph $G=([n],E)$ is {\it left-compressed} if $j_1j_2 \cdots j_r \in E$ implies $i_1i_2 \cdots i_r \in E$ provided $i_p \le j_p$ for every $p, 1\le p\le r$. Equivalently, an $r$-graph $G=([n],E)$ is {\it left-compressed} if $E_{j\setminus i}=\emptyset$ for any $1\le i<j\le n$.
\end{definition}
Recently,  we  proved the following.
\begin{theorem}(Tang et al, 2014) \label{Tang}  Let $m,t,r$ and $p$ be  positive integers satisfying $m={t \choose r}-p$, $0 \leq p \leq t-r+1$, and ${t-r-(p-1) \choose r-2}\geq {t-4 \choose r-4}$  for $p\ge 3$.  Let $G=(V,E)$ be a left-compressed $r$-graph on vertex set $[t]$  with $m$ edges and $|E(G) \triangle E(C_{r,m})|\leq 2$.  Then $\lambda(G)\leq\lambda(C_{r,m})$.
In particular, Conjecture \ref{conjecture} is  true for $r=3$ and ${t \choose 3}-4 \le m \le {t \choose 3}$.
\end{theorem}
The following result implies some Motzkin and Straus'  type theorem holds for   $3$-graphs with   ${t-1 \choose 3} \le m \le {t-1 \choose 3} + {t-2 \choose 2}$ edges.
\begin{theorem}(Peng and Zhao 2013)\label{TheoremPZ18}  Let $m$ and $t$ be positive integers satisfying ${t-1 \choose 3} \le m \le {t-1 \choose 3} + {t-2 \choose 2}$. Let $G$ be a $3$-graph with $m$ edges and  contain a clique of order  $t-1$. Then $\lambda(G) = \lambda\left([t-1]^{(3)}\right)$.
\end{theorem}
Denote $$\lambda_{m}^{r}=\max\{\lambda(G): G {\rm \ is \ an \ } r-{\rm graph\ with \ } m {\rm \ edges }\}.$$
Frankl and F\"uredi's conjecture says that $\lambda_{m}^{r}=\lambda(C_{r,m})$. To verify the truth of Conjecture \ref{conjecture}, it is necessary to show that $\lambda(G)\le \lambda(C_{r,m})$ holds  for every  $r$-graph $G$ with $m$ edges.
Both Lemma \ref{LemmaTal8} and Lemma \ref{LemmaTal9} play a   important role in the exploration of Conjecture \ref{conjecture}.

\begin{lemma}  (Talbot 2002)\label{LemmaTal8}
Let $m,t$ and $r$ be positive integers satisfying  $m\leq{t \choose r}-1$, then there exists a left-compressed $r$-graph $G$ with $m$ edges such that $\lambda(G)=\lambda_{m}^{r}.$
\end{lemma}

\begin{lemma}  (Talbot 2002) \label{LemmaTal9}
Let $m$ and $t$ be positive integers satisfying ${t-1 \choose 3}< m \le {t \choose 3}.$ Let $G=(V,E)$ be a left-compressed $3$-graph with $m$ edges such that $\lambda(G)=\lambda_{m}^{3}$. Let $\vec{x}=(x_{1},x_{2},\ldots ,x_{k})$  be an optimal weighting for $G$ satisfying $x_1 \ge x_2 \ge \cdots \ge x_k >x_{k+1}=\cdots=x_{n}=0$. Then $k\leq t.$
\end{lemma}

For general $r \ge 4$, both Theorems \ref{Addresult} and \ref{Addresult+}   imply $\lambda(G)\leq \lambda(C_{r,m})$ if an $r$-graph $G$ on $[t]$ has $m$ edges and satisfies some nearly structure  of $C_{r,m}$ for some ranges of $m$. We use Theorem \ref{Addresult} in the proof of Theorem \ref{Addresult+}.

\begin{theorem} \label{Addresult} Let $m$, $t$, $r$, $a$ and $i$ be positive integers satisfying  $m={t \choose r}-a$,  where $2i+9\leq a\leq t-r+1$. Let $G=([t],E)$ be a left-compressed $r$-graph  with $m$ edges. If the $r$-tuple with minimum colex ordering in  $G^{c}$ is $(t-r-i+1)(t-r+1)(t-r+3)\ldots t$ then $\lambda(G)\leq\lambda(C_{r,m})$.
\end{theorem}
\begin{theorem} \label{Addresult+} Let $m$, $t$,  $a$ and $r\geq 4$  be positive integers satisfying  $m={t \choose r}-a$,  where $12\leq a\leq t-r+1$. Let $G=([t],E)$ be a left-compressed $r$-graph  with $m$ edges.  If $|E(G)\triangle E(C_{r,m})|\leq 4$ then $\lambda(G)\leq\lambda(C_{r,m})$.
\end{theorem}

The proof of Theorem \ref{Addresult} and  Theorem \ref{Addresult+}  will be given in Section \ref{sect:4}. Next, we state some useful results.

\section{Useful results}
We will impose one additional condition on any optimal weighting ${\vec x}=(x_1, x_2, \ldots, x_n)$ for an $r$-graph $G$:
\begin{eqnarray}
 &&|\{i : x_i > 0 \}|{\rm \ is \ minimal, i.e. \ if}  \ \vec y {\rm \ is \ a \ feasible \ weighting \ for \ } G  {\rm \ satisfying }\nonumber \\
 &&|\{i : y_i > 0 \}| < |\{i : x_i > 0 \}|,  {\rm \  then \ } \lambda (G, {\vec y}) < \lambda(G) \label{conditionb}.
\end{eqnarray}

When the theory of Lagrange multipliers is applied to find the optimum of $\lambda(G, {\vec x})$, subject to $\vec{x}\in S$, notice that $\lambda (E_i, {\vec x})$ corresponds to the partial derivative of  $\lambda(G, \vec x)$ with respect to $x_i$. The following lemma gives some necessary conditions of an optimal weighting for $G$.

\begin{lemma} (Frankl and R\"{o}dl 1984) \label{LemmaTal5} Let $G=(V,E)$ be an $r$-graph on the vertex set $[n]$ and ${\vec x}=(x_1, x_2, \ldots, x_n)$ be an optimal  weighting for $G$ with $k$  ($\le n$) non-zero weights $x_1$, $x_2$, $\ldots$, $x_k$  satisfying condition (\ref{conditionb}). Then for every $\{i, j\} \in [k]^{(2)}$, (a) $\lambda (E_i, {\vec x})=\lambda (E_j, \vec{x})=r\lambda(G)$, (b) there is an edge in $E$ containing both $i$ and $j$.
\end{lemma}

\begin{remark}\label{r1} (a) In Lemma \ref{LemmaTal5}, part (a) implies that
$$x_j\lambda(E_{ij}, {\vec x})+\lambda (E_{i\setminus j}, {\vec x})=x_i\lambda(E_{ij}, {\vec x})+\lambda (E_{j\setminus i}, {\vec x}).$$
In particular, if $G$ is left-compressed, then
$$(x_i-x_j)\lambda(E_{ij}, {\vec x})=\lambda (E_{i\setminus j}, {\vec x})$$
for any $i, j$ satisfying $1\le i<j\le k$ since $E_{j\setminus i}=\emptyset$.

(b) If  $G$ is left-compressed, then for any $i, j$ satisfying $1\le i<j\le k$,
\begin{equation}\label{enbhd}
x_i-x_j={\lambda (E_{i\setminus j}, {\vec x}) \over \lambda(E_{ij}, {\vec x})}
\end{equation}
holds.  If  $G$ is left-compressed and  $E_{i\setminus j}=\emptyset$ for $i, j$ satisfying $1\le i<j\le k$, then $x_i=x_j$.

(c) By (\ref{enbhd}), if  $G$ is left-compressed, then an optimal  weighting  ${\vec x}=(x_1, x_2, \ldots, x_n)$ for $G$  must satisfy
\begin{eqnarray*}\label{conditiona}
x_1 \ge x_2 \ge \cdots \ge x_n \ge 0.
\end{eqnarray*}
\end{remark}




\section{Results on $r$-graphs}\label{sect:4}
In this section we prove Theorems \ref{Addresult} and \ref{Addresult+}.
To prove  Theorem \ref{Addresult}, we need the following lemma.
\begin{lemma}\label{lemmaadd}(Peng et al. Journal of Combinatorial Optimization, accepted) Let $G$ be a left-compressed $r$-graph on the vertex set $[t]$ containing the clique $[t-1]^{(r)}$. Let ${\vec x}=(x_1, x_2, \ldots, x_t)$ be an optimal weighting for $G$. Then
\begin{eqnarray*}
 x_1\le x_{t-1}+x_t \le 2x_{t-1}.
\end{eqnarray*}
\end{lemma}
\noindent{\em Proof of Theorem \ref{Addresult}.}
Recall that the $r$-tuple with minimum colex ordering in  $G^{c}$ is $(t-r-i+1)(t-r+1)(t-r+3)\ldots t$ and $m={t \choose r}-a$, where $2i+9\leq a\leq t-r+1$. In view of Fig.1, we have
$$G=[t]^{(r)}\backslash \left\{\bigcup_{j=t-r+2-a+i}^{t-r+1}j(t-r+2)\ldots t,\bigcup_{j=1}^{i}(t-r-j+1)(t-r+1)(t-r+3)\ldots t\right\}$$
and
$$C_{r,m}=[t]^{(r)}\backslash \left\{\bigcup_{j=t-r+2-a}^{t-r+1}j(t-r+2)\ldots t\right\}.$$
Let $\vec{x}=(x_{1},x_{2},\ldots ,x_{t})$  be an optimal weighting for $G$ satisfying $x_1 \ge x_2 \ge \cdots \ge x_t \ge 0$.
If $x_t=0$ then $\lambda(G)\leq\lambda([t-1]^{(r)})$ since $\vec{x}$ has at most $t-1$ positive weights. Since  $m={t \choose r}-a\geq{t-1 \choose r}$ for $a\leq t-r+1$, then $[t-1]^{(r)}\subseteq C_{r,m}$. Hence $\lambda(C_{r,m})\geq  \lambda([t-1]^{(r)})\geq \lambda (G)$ by Fact \ref{mono}. Therefore Theorem\ref{Addresult} holds. So we assume $x_t>0$ next.

First we point out that
\begin{eqnarray}\label{ineqadd}
4\lambda(E_{1(t-r-i+1)},\vec{x})\geq\lambda(E_{(t-r+1)(t-r+2)},\vec{x}).
\end{eqnarray}

To show this, note that
\begin{eqnarray*}\lambda(E_{1(t-r-i+1)},\vec{x})&=&\sum\limits_{j_1,j_2,\ldots, j_{r-2}\in [t]\setminus \{1,t-r+1-i,t-r+1,t-r+2\}}x_{j_1}x_{j_2}\ldots x_{j_{r-2}}\\
&+&x_{t-r+1}\sum\limits_{j_1,j_2,\ldots, j_{r-3}\in [t]\setminus \{1,t-r+1-i,t-r+1, t-r+2\}}x_{j_1}x_{j_2}\ldots x_{j_{r-3}}\\
&+&x_{t-r+2}\sum\limits_{j_1,j_2,\ldots, j_{r-3}\in [t]\setminus \{1,t-r+1-i,t-r+1, t-r+2\}}x_{j_1}x_{j_2}\ldots x_{j_{r-3}}\\
&+&x_{t-r+1}x_{t-r+2}\sum\limits_{j_1,j_2,\ldots, j_{r-4}\in [t]\setminus \{1,t-r+1-i,t-r+1,t-r+2\}}x_{j_1}x_{j_2}\ldots x_{j_{r-4}}
\end{eqnarray*}
and
\begin{eqnarray*}\lambda(E_{(t-r+1)(t-r+2)},\vec{x})&=&\sum\limits_{j_1,j_2,\ldots, j_{r-2}\in [t]\setminus \{1,t-r-i+1,t-r+1,t-r+2\}}x_{j_1}x_{j_2}\ldots x_{j_{r-2}}\\
&+&x_{1}\sum\limits_{j_1,j_2,\ldots, j_{r-3}\in [t]\setminus \{1,t-r-i+1,t-r+1,t-r+2\}}x_{j_1}x_{j_2}\ldots x_{j_{r-3}}\\
&+&x_{t-r-i+1}\sum\limits_{j_1,j_2,\ldots, j_{r-3}\in [t]\setminus \{1,t-r-i+1,t-r+1,t-r+2\}}x_{j_1}x_{j_2}\ldots x_{j_{r-3}}\\
&+&x_{1}x_{t-r-i+1}\sum\limits_{j_1,j_2,\ldots, j_{r-4}\in [t]\setminus \{1,t-r-i+1,t-r+1,t-r+2\}}x_{j_1}x_{j_2}\ldots x_{j_{r-4}}-x_{t-r+3}\ldots x_{t-1}x_t.
\end{eqnarray*}
Since $G$ contains $[t-1]^{(r)}$, we have $x_1\leq x_{t-1}+x_t\leq 2 x_{t-r+1}$ and $x_{t-r+1}\leq x_1\leq x_{t-1}+x_t\leq 2x_{t-r+2}$ by Lemma \ref{lemmaadd}. So
\begin{eqnarray*}\lambda(E_{(t-r+1)(t-r+2)},\vec{x})&\leq &\sum\limits_{j_1,j_2,\ldots, j_{r-2}\in [t]\setminus \{1,t-r-i+1,t-r+1,t-r+2\}}x_{j_1}x_{j_2}\ldots x_{j_{r-2}}\\
&+& 2x_{t-r+1}\sum\limits_{j_1,j_2,\ldots, j_{r-3}\in [t]\setminus \{1,t-r-i+1,t-r+1,t-r+2\}}x_{j_1}x_{j_2}\ldots x_{j_{r-3}}\\
&+&2x_{t-r+2}\sum\limits_{j_1,j_2,\ldots, j_{r-3}\in [t]\setminus \{1,t-r-i+1,t-r+1,t-r+2\}}x_{j_1}x_{j_2}\ldots x_{j_{r-3}}\\
&+&4x_{t-r+1}x_{t-r+2}\sum\limits_{j_1,j_2,\ldots, j_{r-4}\in [t]\setminus \{1,t-r-i+1,t-r+1,t-r+2\}}x_{j_1}x_{j_2}\ldots x_{j_{r-4}}-x_{t-r+3}\ldots x_{t-1}x_t\\
&\leq&4\lambda(E_{1(t-r-i+1)},\vec{x}).
\end{eqnarray*}
So (\ref{ineqadd}) holds. Using (\ref{ineqadd}), we prove that
\begin{eqnarray}\label{xy}
x_{t-r+1}-x_{t-r+2}\geq x_{1}-x_{t-r-i+1}.
\end{eqnarray}
Recall that $x_1 \ge x_2 \ge \cdots \ge x_t > 0$.
By Remark \ref{r1} (b), we have
\begin{eqnarray*}\label{eqII}x_{1}&=&x_{t-r-i+1}+\frac{\lambda(E_{1\setminus (t-r-i+1)},\vec{x})}{\lambda(E_{1(t-r-i+1)},\vec{x})}\nonumber \\
&=&x_{t-r-i+1}+\frac{(x_{t-r+1}+x_{t-r+2})x_{t-r+3}\ldots x_{t}}{\lambda(E_{1(t-r-i+1)},\vec{x})},
\end{eqnarray*}
and
\begin{eqnarray*}x_{t-r+1}&=&x_{t-r+2}+\frac{\lambda\left(E_{(t-r+1)\setminus (t-r+2)},\vec{x}\right)}{\lambda(E_{(t-r+1)(t-r+2)},\vec{x})}\nonumber \\
&=&x_{t-r+1}+\frac{(x_{t-r-a+i+2}+\cdots +x_{t-r-i})x_{t-r+3}\ldots x_{t}}{\lambda(E_{(t-r+1)(t-r+3)},\vec{x})}.
\end{eqnarray*}
Recall that $a\geq 2i+9$, $x_{t-r-a+i+2}  \ge \cdots \ge x_{t-r-i}\geq x_{t-r+1}\geq x_{t-r+2}$ and $\lambda(E_{(t-r+1)(t-r+2)},\vec{x})\leq 4\lambda(E_{1(t-r-i+1)},\vec{x})$.
We have
\begin{eqnarray*}x_{t-r+1}-x_{t-r+2}&=&x_{t-r+1}+\frac{(x_{t-r-a+i+2}+\cdots +x_{t-r-i})x_{t-r+3}\ldots x_{t}}{\lambda(E_{(t-r+1)(t-r+3)},\vec{x})}\\
&\geq &\frac{4(x_{t-r+1}+x_{t-r+2})x_{t-r+3}\ldots x_{t}}{4\lambda(E_{1(t-r-i+1)},\vec{x})}\\
&=&x_{1}-x_{t-r-i+1}.
\end{eqnarray*}
Hence (\ref{xy}) holds.

Using (\ref{xy}), we prove that $\lambda(G) \leq \lambda(C_{r,m})$. We have $x_{1}=x_{2}=\cdots=x_{t-r-a+i+1}$ and $x_{t-r-i+1}=\cdots=x_{t-r}$ by Remark \ref{r1} (b). Hence \begin{eqnarray*}\lambda(C_{r,m},\vec{x})-\lambda(G,\vec{x})&=&i(x_{t-r-i+1}x_{t-r+1}x_{t-r+3}\ldots x_{t}-x_{1}x_{t-r+2}\ldots x_{t})\\
&=&i(x_{t-r-i+1}x_{t-r+1}-x_{1}x_{t-r+2})x_{t-r+3}\ldots x_{t}.
\end{eqnarray*}
\begin{eqnarray*}\label{eq37}
\lambda(C_{r,m},\vec{x})-\lambda(G,\vec{x})&=&i(x_{t-r-i+1}x_{t-r+1}-x_{1}x_{t-r+2})x_{t-r+3}\ldots x_{t}\nonumber \\
&=& i[x_{t-r-i+1}(x_{t-r+1}+x_{t-r+2}-x_{t-r+2})-x_{1}x_{t-r+2}]x_{t-r+3}\ldots x_{t} \nonumber \\
&\geq& i[ x_{t-r-i+1}(x_{t-r+2}+x_{1}-x_{t-r-i+1})-x_{1}x_{t-r+2}]x_{t-r+3}\ldots x_{t}\nonumber \\
&=&i(x_{t-r-i+1}-x_{t-r+2})(x_{1}-x_{t-r-i+1})x_{t-r+3}\ldots x_{t}\nonumber \\
&\geq& 0.
\end{eqnarray*}
Therefore $\lambda(C_{r,m})\geq \lambda(C_{r,m},\vec{x})\geq \lambda(G, \vec{x})=\lambda(G)$. This completes the proof of Theorem \ref{Addresult}. \qed

To prove  Theorem \ref{Addresult+}, we need the following lemma.
\begin{lemma}\label{lemmaadd+} Let $m$, $t$, $a$ and  $r\geq 4$  be positive integers satisfying  $m={t \choose r}-a$,  where $12\leq a\leq t-r+1$. Let $$G=[t]^{(r)}\backslash \left\{\bigcup_{j=t-r-a+4}^{t-r+1}j(t-r+2)\ldots t, (t-r)(t-r+1)(t-r+3)\ldots t, (t-r)(t-r+1)(t-r+2)(t-r+4)\ldots t\right\}.$$ Then $\lambda(G)\leq \lambda(C_{r,m}).$
\end{lemma}
\noindent{\em Proof.} The proof of this lemma is similar to the proof of Theorem \ref{Addresult}. Let $\vec{x}=(x_{1},x_{2},\ldots ,x_{t})$  be an optimal weighting for $G$ satisfying $x_1 \ge x_2 \ge \cdots \ge x_t \ge 0$. Similar as we did in the proof of Theorem \ref{Addresult}, we have if $x_t=0$ then $\lambda(G)\leq \lambda (C_{r,m})$. So we  assume that $x_t>0$ next.

First we point out that
\begin{eqnarray}\label{ineqadd+}
4\lambda(E_{(t-r-a+3)(t-r)},\vec{x})\geq\lambda(E_{(t-r+1)(t-r+3)},\vec{x}).
\end{eqnarray}

To show this, note that
\begin{eqnarray*}\lambda(E_{(t-r-a+3)(t-r)},\vec{x})&=&\sum\limits_{j_1,j_2,\ldots, j_{r-2}\in [t]\setminus \{t-r-a+3,t-r,t-r+1,t-r+3\}}x_{j_1}x_{j_2}\ldots x_{j_{r-2}}\\
&+&x_{t-r+1}\sum\limits_{j_1,j_2,\ldots, j_{r-3}\in [t]\setminus \{t-r-a+3,t-r,t-r+1,t-r+3\}}x_{j_1}x_{j_2}\ldots x_{j_{r-3}}\\
&+&x_{t-r+3}\sum\limits_{j_1,j_2,\ldots, j_{r-3}\in [t]\setminus \{t-r-a+3,t-r,t-r+1,t-r+3\}}x_{j_1}x_{j_2}\ldots x_{j_{r-3}}\\
&+&x_{t-r+1}x_{t-r+3}\sum\limits_{j_1,j_2,\ldots, j_{r-4}\in [t]\setminus\{t-r-a+3,t-r,t-r+1,t-r+3\}}x_{j_1}x_{j_2}\ldots x_{j_{r-4}}
\end{eqnarray*}
and
\begin{eqnarray*}\lambda(E_{(t-r+1)(t-r+3)},\vec{x})&=&\sum\limits_{j_1,j_2,\ldots, j_{r-2}\in [t]\setminus \{t-r-a+3,t-r,t-r+1,t-r+3\}}x_{j_1}x_{j_2}\ldots x_{j_{r-2}}\\
&+&x_{t-r-a+3}\sum\limits_{j_1,j_2,\ldots, j_{r-3}\in [t]\setminus \{t-r-a+3,t-r,t-r+1,t-r+3\}}x_{j_1}x_{j_2}\ldots x_{j_{r-3}}\\
&+&x_{t-r}\sum\limits_{j_1,j_2,\ldots, j_{r-3}\in [t]\setminus \{t-r-a+3,t-r,t-r+1,t-r+3\}}x_{j_1}x_{j_2}\ldots x_{j_{r-3}}\\
&+&x_{t-r-a+3}x_{t-r}\sum\limits_{j_1,j_2,\ldots, j_{r-4}\in [t]\setminus \{t-r-a+3,t-r,t-r+1,t-r+3\}}x_{j_1}x_{j_2}\ldots x_{j_{r-4}}\\
&-&x_{t-r+2}x_{t-r+4}\ldots x_{t-1}x_t-x_{t-r}x_{t-r+4}\ldots x_{t-1}x_t.
\end{eqnarray*}
Note that $r\geq 4$. In view of Fig.1, we have $G$ contains $[t-1]^{(r)}.$ So $x_{t-r-a+3}\leq x_1\leq x_{t-1}+x_t\leq 2 x_{t-r+1}$ and $x_{t-r}\leq x_1\leq x_{t-1}+x_t\leq 2x_{t-r+3}$ by Lemma \ref{lemmaadd}. Hence
\begin{eqnarray*}\lambda(E_{(t-r+1)(t-r+3)},\vec{x})&\leq &\sum\limits_{j_1,j_2,\ldots, j_{r-2}\in [t]\setminus \{1,t-r-i+1,t-r+1,t-r+2\}}x_{j_1}x_{j_2}\ldots x_{j_{r-2}}\\
&+& 2x_{t-r+1}\sum\limits_{j_1,j_2,\ldots, j_{r-3}\in [t]\setminus \{1,t-r-i+1,t-r+1,t-r+2\}}x_{j_1}x_{j_2}\ldots x_{j_{r-3}}\\
&+&2x_{t-r+3}\sum\limits_{j_1,j_2,\ldots, j_{r-3}\in [t]\setminus \{1,t-r-i+1,t-r+1,t-r+2\}}x_{j_1}x_{j_2}\ldots x_{j_{r-3}}\\
&+&4x_{t-r+1}x_{t-r+3}\sum\limits_{j_1,j_2,\ldots, j_{r-4}\in [t]\setminus \{1,t-r-i+1,t-r+1,t-r+2\}}x_{j_1}x_{j_2}\ldots x_{j_{r-4}}\\
&-&x_{t-r+2}x_{t-r+4}\ldots x_{t-1}x_t-x_{t-r}x_{t-r+4}\ldots x_{t-1}x_t\\
&\leq&4\lambda(E_{(t-r-a+3)(t-r)},\vec{x}).
\end{eqnarray*}
So (\ref{ineqadd+}) holds. Using (\ref{ineqadd+}), we prove that
\begin{eqnarray}\label{xyz}
x_{t-r+1}-x_{t-r+3}\geq x_{t-r-a+3}-x_{t-r}.
\end{eqnarray}
Recall that $x_1 \ge x_2 \ge \cdots \ge x_t > 0$.  By Remark \ref{r1} (b), we have
\begin{eqnarray*}\label{eqII}x_{t-r-a+3}&=&x_{t-r}+\frac{\lambda(E_{(t-r-a+3)\setminus (t-r)},\vec{x})}{\lambda(E_{(t-r-a+3)(t-r)},\vec{x})}\nonumber \\
&=&x_{t-r-i+1}+\frac{(x_{t-r+1}+x_{t-r+2})x_{t-r+3}\ldots x_{t}}{\lambda(E_{(t-r-a+3)(t-r)},\vec{x})},
\end{eqnarray*}
and
\begin{eqnarray*}x_{t-r+1}&=&x_{t-r+3}+\frac{\lambda\left(E_{(t-r+1)\setminus (t-r+3)},\vec{x}\right)}{\lambda(E_{(t-r+1)(t-r+3)},\vec{x})}\nonumber \\
&=&x_{t-r+1}+\frac{(x_{t-r-a+4}+\cdots +x_{t-r-1})x_{t-r+2}x_{t-r+4}\ldots x_{t}}{\lambda(E_{(t-r+1)(t-r+3)},\vec{x})}.
\end{eqnarray*}
Recall that $a\geq 12$, $x_{t-r-a+4}\geq\cdots \geq x_{t-r-1}\geq x_{t-r+1}\geq x_{t-r+2}$ and $4\lambda(E_{(t-r-a+3)(t-r)},\vec{x})\geq\lambda(E_{(t-r+1)(t-r+3)},\vec{x}).$
We have
\begin{eqnarray*}x_{t-r+1}-x_{t-r+3}&=&x_{t-r+1}+\frac{(x_{t-r-a+4}+\cdots +x_{t-r-1})x_{t-r+2}x_{t-r+4}\ldots x_{t}}{\lambda(E_{(t-r+1)(t-r+3)},\vec{x})}\\
&\geq& \frac{4(x_{t-r+1}+x_{t-r+2})x_{t-r+3}\ldots x_{t}}{4\lambda(E_{(t-r-a+3)(t-r)},\vec{x})}\\
&\geq& x_{t-r-a+3}-x_{t-r}.
\end{eqnarray*}
Hence (\ref{xyz}) holds. Using (\ref{xyz}), we prove that $\lambda(G)\leq \lambda(C_{r,m}).$
\begin{eqnarray*}\label{eq37}
\lambda(C_{r,m},\vec{x})-\lambda(G,\vec{x})&=&(x_{t-r}x_{t-r+1}-x_{t-r-a+3}x_{t-r+3})x_{t-r+2}x_{t-r+4}\ldots x_{t}\nonumber \\
&=& [x_{t-r}(x_{t-r+1}+x_{t-r+3}-x_{t-r+3})-x_{t-r-a+3}x_{t-r+3}]x_{t-r+2}x_{t-r+4}\ldots x_{t} \nonumber \\
&\geq& [x_{t-r}(x_{t-r+3}+x_{t-r-a+3}-x_{t-r})-x_{t-r-a+3}x_{t-r+3}]x_{t-r+2}x_{t-r+4}\ldots x_{t}\nonumber \\
&=&(x_{t-r}-x_{t-r+3})(x_{t-r-a+3}-x_{t-r})x_{t-r+2}x_{t-r+4}\ldots x_{t}\nonumber \\
&\geq& 0.
\end{eqnarray*}
Therefore $\lambda(C_{r,m})\geq \lambda(C_{r,m},\vec{x})\geq \lambda(G, \vec{x})=\lambda(G)$. This completes the proof of Lemma \ref{lemmaadd+}.\qed

Now we are ready to prove Theorem \ref{Addresult+}.

\noindent{\em Proof of Theorem \ref{Addresult+}.} Since $|E(G) \Delta E(C_{r,m})| \le 4$ and  $G$ is left-compressed on $[t]$, in view of Fig.1, there are only the following possible cases for $E(G)$ except for $E(G)=C_{r,m}$:

Case 1. $$E(G)=[t]^{(r)}\backslash \left\{\bigcup_{j=t-r-a+3}^{t-r+1}j(t-r+2)\ldots t, (t-r)(t-r+1)(t-r+3)\ldots t\right\};$$

Case 2. $$E(G)=[t]^{(r)}\backslash \left\{\bigcup_{j=t-r-a+4}^{t-r+1}j(t-r+2)\ldots t, (t-r)(t-r+1)(t-r+3)\ldots t, (t-r-1)(t-r+1)(t-r+2)(t-r+4)\ldots t\right\};$$

Case 3. $$E(G)=[t]^{(r)}\backslash \left\{\bigcup_{j=t-r-a+4}^{t-r+1}j(t-r+2)\ldots t, (t-r)(t-r+1)(t-r+3)\ldots t, (t-r)(t-r+1)(t-r+2)(t-r+4)\ldots t\right\}.$$

Cases 1 and 2 correspond to $i=1$ and $i=2$ in Theorem \ref{Addresult}, respectively. So if Cases 1 and 2 happen, then $\lambda(G)\leq\lambda(C_{r,m})$ by Theorem \ref{Addresult}. If Case 3   happens, then by Lemma \ref{lemmaadd+}, $\lambda(G)\leq\lambda(C_{3,m})$.  The proof of Theorem \ref{Addresult+} is completed.
\qed

\bibliographystyle{spmpsci}      
\bibliographystyle{unsrt}

\end{document}